\documentclass[a4paper,12pt,leqno]{article}
\usepackage{latexsym}
\usepackage[all]{xy}

\usepackage{amssymb} 
\usepackage{amsmath} 

\usepackage{tikz}
\usetikzlibrary{arrows,cd}  
\usepackage{amsthm}

\def\P{{\mathbf{P}}}
\def\Z{{\mathbb{Z}}}
\def\Q{{\mathbb{Q}}}
\def\K{{\mathbb{K}}}
\def\CC{{\mathbb{C}}}
\def\R{{\mathbb{R}}}
\def\A{{\mathcal{A}}}
\def\B{{\mathcal{B}}}

\def\O{{\mathcal{O}}}

\DeclareMathOperator{\codim}{codim}

\DeclareMathOperator{\Der}{Der}

\DeclareMathOperator{\pd}{pd}

\DeclareMathOperator{\reg}{reg}

\numberwithin{equation}{section}

\newcommand{\owari}{\hfill$\square$}

\newtheorem{theorem}{Theorem}[section]
\newtheorem{prop}[theorem]{Proposition}
\newtheorem{cor}[theorem]{Corollary}

\newtheorem{define}[theorem]{Definition}

\newtheorem{conj}[theorem]{Conjecture}
\theoremstyle{remark}
\newtheorem{rem}[theorem]{Remark}
\newtheorem{example}[theorem]{Example}

\title{Solomon-Terao polynomials and Castelnouvo-Mumford regularity of hyperplane arrangements}
\author{Takuro Abe}
\date{\today}

\begin{document}

\maketitle

\begin{abstract}
The Solomon-Terao bi-polynomial was introduced by Solomon and Terao which degenerates to the characteristic polynomial of hyperplane arrangements. Also, it was proved recently that the other specialization of the Solomon-Terao bi-polynomial, we call the Solomon-Terao polynomial, 
coincides with the Poinrar\'{e} polynomial of the regular nilpotent Hessenberg variety when the arrangement and the variety comes from the same lower ideal in the positive system. Moreover, there are recent developments with superspace coinvariants and Fields conjecture, thus these polynomials are becoming more and more important. However, the research of them has been very hard, and even the top degree of the Solomon-Terao polynomial has not yet been known, which we solve in this article, by using the Castelnouvo-Mumford regularity of the logarithmic derivation modules.
 \end{abstract}

\section{Introduction}
Let $\K$ be a field, $V=\K^\ell$ and $\A$ be an irreducible arrangement in $V$. Let $V^*=
\langle x_1,\ldots,x_\ell\rangle_\K$ and let $S=\mbox{Sym}^*(V^*)=\K[x_1,\ldots,x_\ell]$ be the coordinate ring of $V$. Let $\Der S:=\bigoplus_{i=1}^\ell S \partial_{x_i}
$ be the $S$-graded free module of $S$-derivations of rank $\ell$, where let us agree that $\deg \partial_{x_i}=0$. 
For each $H \in \A$ let us fix a linear form $\alpha_H \in V^*$ such that 
$\ker \alpha_H=H$, and let $Q(\A):=\prod_{H \in \A} \alpha_H$. Now let us introduce two fundamental invariants of hyperplane arrangements.

\begin{define}

(1)\,\,
Let 
$$
L(\A):=\{\bigcap_{H \in \B} H \mid \B \subset \A\}
$$
be the \textbf{intersection lattice} of $\A$, and let $\mu:L(\A) \rightarrow 
\Z$ be the \textbf{M\"{o}bius function} defined by 
$\mu(V)=1$ and by 
$$
\mu(X):=-\sum_{X \subsetneq Y \subset V,\ Y \in L(\A)} \mu(Y)
$$
for $X \neq V$. Let $$
L_k(\A):=\{X \in L(\A) \mid \codim_V X=k\}.
$$ 

(2)\,\,
The \textbf{characteristic polynomial} $\chi(\A;t)$ of $\A$ is defined by 
$$
\chi(\A;t):=\sum_{X \in L(\A)} \mu(X) t^{\dim X},
$$
and the \textbf{Poincar\`{e} polynomial} $\pi(\A;t)$ of $\A$ is defined by 
$$
\pi(\A;t):=\sum_{X \in L(\A)} \mu(X) (-t)^{\codim X}.
$$
\label{lattice}
\end{define}

\begin{define}
(1)\,\,
For $ 0\le p \le \ell$, let 
$$
D^p(\A):=\{ \theta \in \wedge^p \Der S \mid 
\theta(\alpha_H,f_2,\ldots,f_p) \in S \alpha_H\ (\forall H \in \A,\ \forall f_i \in S)\}
$$
be the \textbf{logarithmic derivation module of order $p$}.

(2)\,\,
We say that $\A$ is \textbf{free} with exponents $\exp(\A)=(d_1,\ldots,d_\ell)$ if 
$$
D(\A) \simeq \bigoplus_{i=1}^\ell S[-d_i].
$$
\label{logmodule}
\end{define}

Clearly Definition \ref{lattice} gives us combinatorial invariants, and Definition \ref{logmodule} 
 algebraic invariants. In fact, they are related to each other, and also to geometry. For example, when $\K=\CC$, then 
$$
\mbox{Poin}(V \setminus \cup_{H \in \A}H;t)= \pi(\A;t)
$$
and when $\A$ is free with exponents $(d_1,\ldots,d_\ell)$, then 
$$
\pi(\A;t)=\prod_{i=1}^\ell (1+d_it).
$$
The former is shown by Orlik-Solomon in \cite{OS}, and the latter by Terao in \cite{T2}. In this article we study a different algebraic invariants coming from the higher order logarithmic derivations. To explain it, let us recall some results from the wonderful paper by Solomon and Terao \cite{ST}. 
Let 
$$
\Psi(\A;x,t):
=\sum_{p=0}^\ell \mbox{Hilb}(D^p(\A);x)(t(x-1)-1)^p 
$$
be a series. Then in \cite{ST}, the following surprising results were proved.

\begin{theorem}[Theorem 1.2 and 
Proposition 5.3, \cite{ST}]
The series $\Psi(\A;x,t)$ is in fact a polynomial in $ \Q[x,t]$. Moreover, 
$$
(-1)^\ell\Psi(\A;1,t)=\chi(\A;t).
$$
\label{ST}
\end{theorem}

So from higher order logarithmic derivation modules we can recover the characteristic polynomial. However, in the same level, the surprising result is that $\Psi(\A;x,t)$ is a bi-polynomial in $x$ and $t$, and what Theorem \ref{ST} says is that $\Psi(\A;x.t)$ has a great specialization of $x=1$. Then it is natural to consider the specialization of $t$. To see it let us summarize and name what appeared here. 

\begin{define}[\cite{AMMN}]

(1)\,\,
$\Psi(\A;x,t)$ is called the \textbf{Solomon-Terao bi-polynomial}. 
of $\A$.

(2)\,\,
$$
ST(\A;x):=\Psi(\A;x,-1)
= \sum_{p=0}^\ell \mbox{Hilb}(D^p(\A);x)(-x)^p \in \Q[x]
$$
is called the \textbf{Solomon-Terao polynomial} of $\A$.
\label{STpoly}
\end{define}

Also surprisingly, there is a very interesting relation between $ST(\A;x)$ and geometry.
For several definitions and terminologies, see \cite{AHMMS}. 

\begin{theorem}[Theorem 1.3 and Corollary 7.3, \cite{AHMMS}]
Assume that $\A=\A_I$ is an ideal subarrangement of the Weyl arrangement corresponding to the lower ideal $I$ of the positive system $\Phi^+$, see \cite{AHMMS} for details. Then 
$ST(\A_I;x)$ coincides with the topological Poincar\`{e} polynomial of the regular nilpotent 
Hessenberg variety defined by $I$.
\label{Hessenberg}
\end{theorem}

Hence the Solomon-Terao bi-polynomial has two very nice specializations which have geometric interpretations. Also, recently, there have been found applications of Solomon-Terao polynomials and algebras (introduced in \cite{AMMN}, see \S3 for its generalization to multiarrangements) to superspace coinvariants in \cite{ACKMR}, which solved the Fields conjecture in \cite{MRW}. So it is very natural to ask more general meaning of the Solomon-Terao polynomial when the arrangement is not necessarily an ideal arrangement. However, we know almost nothing on $ST(\A;x)$ in 
general. For example, we do not know the top degree term of this polynomial, and its degree in fact when $\ell \ge 3$! What we know is that, when $\A$ is free with exponents $(d_1,\ldots,d_\ell)$, i.e., when 
$$
D(\A) \simeq \oplus_{i=1}^\ell S[-d_i],
$$
then 
$$
\Psi(\A;x,t)=\prod_{i=1}^\ell (1+x+\cdots+x^{d_i-1}-tx^{d_i})
$$
and thus 
$$
ST(\A;x)=\Psi(\A;x,-1)=\prod_{i=1}^\ell (1+x+\cdots+x^{d_i}).
$$
Also a recent result in \cite{A12} gives us a way to compute $ST(\A;x)$ inductively when $\A$ and the restriction of $\A$ onto some hyperplane are both free. So to study these polynomials, the research on the fundamentals on these polynomials are very important.

Based on many computations and formulae of special cases in \cite{AMMN} and \cite{AHMMS}, the following is conjectured:

\begin{conj}[\cite{AMMN}, Conjecture 5.12]
The polynomial $\Psi(\A;x,-1)$ is a monic polynomial of degree $|\A|$ when $\A$ is tame.
\label{conj1}
\end{conj}

Here we say that $\A$ 
is \textbf{tame} if $\pd_S D^{\ell-p}(\A)^* \le p$ for all $0 \le p \le \ell$. See \cite{A14} for 
the recent developments of tame arrangement theory. The first main result in this article is to settle this conjecture affirmatively:

\begin{theorem}
Conjecture \ref{conj1} is true.
\label{main1}
\end{theorem}

Since all $3$-arrangements are tame, we have 

\begin{cor}
Let $\ell=3$. Then the top degree term of $ST(\A;x)=\Psi(\A;x,-1)$ is 
$ x^{|\A|}$.
\label{3topdegree}
\end{cor}

Also we can get data for the second largest coefficients of the Solomon-Terao polynomials.

\begin{theorem}
Let $\A$ be tame and irreducible. Then 
$$
ST(\A;x)=\Psi(\A;x,-1)=1+\ell x+\cdots+(\ell+a)x^{n-1}+x^n
$$
with $n:=|\A|,\ a \ge 0$. Here  $a$ is the number of  
relations of degree $|\A|$ among a minimal set of generators for $D^{\ell-1}(\A)$.
\label{main2}
\end{theorem}

In fact we prove these results in the class of multiarrangements $(\A,m)$, for the definition see \S2. Since $\A$ is a special case of multiarrangements, we can show them in a wider class.

The main tool to show Theorem \ref{main1} is the famous (Castelnouvo-Mumford) regularity. We study that of $D^p(\A,m)$ and $\Omega^p(\A,m)$ to show
Theorem \ref{main1} (see Definition \ref{multilog} for their definitions). One of the main tools to show it is the following 
bound for the regularity, which we can show by applying the similar argument as in \cite{MS}.
Note that this proof is basically same as that by Morihiko Saito in \cite{MS} for $D(\A)$, so the following is a generalization to higher order multiarrangements.

\begin{theorem}
Assume that $(\A,m)$ is essential and let $|m|=n$. Then 
$$
\reg^p(\A,m):=\reg(D^p(\A,m)) \le |m|-\ell+p
$$
for $0 \le p \le \ell$. Also,
$$
\reg(\Omega^p(\A,m)) \le -p
$$
for $0 \le p \le \ell$. 
\label{regD}
\end{theorem}

\begin{rem}
In \cite{MS} the degree of $\partial_{x_i}$ is $-1$, but here it is $0$. So the statement is the same though the expressions are different. 

\end{rem}

\textbf{Acknowledgements}.
The author 
is partially supported by JSPS KAKENHI Grant Numbers JP23K17298 and JP25H00399. We are grateful to Lukas K\"{u}hne for his 
computation of the logarithmic differential forms in Example \ref{ex2}.

\section{Preliminaries}

Let us summarize several results used for the proof of the main theorems. First let us define the dual of $D^p(\A)$.

\begin{define}
Let $\Omega^1_V:=\oplus_{i=1}^\ell Sdx_i $ and let $\Omega^p_V:=\wedge^p \Omega^1_V$. Then 
$$
\Omega^p(\A):=\{ \omega \in 
\frac{1}{Q(\A)} \Omega^p_V \mid Q(\A) \omega\wedge d\alpha_H \in \alpha_H \Omega^{p+1}_V\ (\forall 
H \in \A)\}
$$
is the module of \textbf{logarithmic differential $p$-forms}.
\label{logdiff}
\end{define}

It is known that 
$$
D^p(\A)^* \simeq \Omega^p(\A),\ 
D^p(\A) \simeq \Omega^p(\A)^*.
$$
So both are reflexive module, and $\pd_S D^p(\A) \le \ell-2,\ 
\pd_S \Omega^p(\A) \le \ell-2$ by the reflexivity and the Auslander-Buchsbaum formula. Next recall the definition of the Castelnouvo-Mumford regularity.

\begin{define}
Let $M$ be a finitely generated $S$-graded module. Let 
$$
0 \rightarrow F_m \rightarrow F_{m-1} \rightarrow 
\cdots F_0 \rightarrow M \rightarrow 0
$$
be a graded minimal free resolution of $M$ with 
$$
F_i:=\oplus_{j=1}^{s_i} S[-d_{ij}].
$$
Then the \textbf{(Castelnouvo-Mumford) regularity $\reg M$} is 
$$
\reg M:=\max\{ d_{ij}-i \mid 0 \le i \le m,\ 1 \le j \le s_i\}.
$$
\label{reg}
\end{define}

On the regularity, the furst result is due to Schenck in \cite{Sc} for $\reg (D(\A))$ when $\ell=3$, which are now generalized by Morihiko Saito and 
Bath as follows:

\begin{prop}[\cite{MS}, Proposition 1.3]
Let $\reg^p(\A):=\mbox{reg}(D^p(\A))$ denote the 
Castelnouvo-Mumford regularity of $D^p(\A)$. Then 
$\reg^1(\A) \le |\A|-\ell+1$.
\label{MS}
\end{prop}

\begin{theorem}[\cite{B}]
$$
\reg(\Omega^p(\A)) \le -p
$$
for $0 \le p \le \ell$.
\label{Bath}
\end{theorem}

\begin{rem}
In \cite{B}, $\deg d{x_i}=1$, but in this article 
$\deg d{x_i}=0$. So the formulations in \cite{B} and Theorem \ref{Bath} are different.
\end{rem}

The following approximation-type result due to Derksen and Sidman plays the key role in the study of regularity of logarithmic modules as in \cite{MS}.

\begin{theorem}[\cite{DS}]
Let $M \subset \Der S$ be a graded $S$-module. Assume that there are 
graded $S$-modules $M_1,\ldots,M_s \subset \Der S$ and ideals $I_1,\ldots,I_s \subset S$ such that 
$$
I_i M_i \subset M \subset M_i
$$
for all $i$, and $I_1+\cdots+I_s=(x_1,\ldots,x_\ell)$. Then $\reg (M) \le r$ if $\reg (M_i) \le r-1$ for all $i$ and some $r \ge 2$.
\label{DS}
\end{theorem}

Now let us introduce a multiarrangement. We say that $(\A,m)$ is a 
\textbf{multiarrangement} if $m$ is a function from $\A$ to $\Z_{>0}$. Let $Q(\A,m):=
\prod_{H \in \A} \alpha_H^{m(H)}$ and $|m|:=
\deg Q(\A,m)$. 
Then we can define their logarithmic modules as follows.

\begin{define}
The \textbf{logarithmic derivation module $D^p(\A,m)$ 
of order $p$} of $(\A,m)$ is defined by 
$$
D^p(\A,m):=\{\theta \in \wedge^p \Der S \mid 
\theta(\alpha_H,f_2,\ldots,f_p) \in S\alpha_H^m(H)\ (\forall 
H \in \A,\ f_i \in S)\}.
$$
The \textbf{logarithmic differential $p$-forms} $ \Omega^p(\A,m)$ of $(\A,m)$ is defined by 
$$
\Omega^p(\A,m):=\{\omega \in \frac{1}{Q(\A,m)}\Omega^p_V \mid 
Q(\A,m)\omega \wedge d\alpha_H \in S\alpha_H^{m(H)}\ (\forall 
H \in \A,\ f_i \in S)\}.
$$
\label{multilog}
\end{define}

\begin{prop}[\cite{Z}]
$D^p(\A,m)$ and $\Omega^p(A,m)$ are $S$-dual modules, so they are reflexive. In particular, their projective dimensions are at most $\ell-2$.
\label{dualmulti}
\end{prop}

\begin{define}
(1)\,\,
We say that $(\A,m)$ is \textbf{free} with $\exp(\A,m)=(d_1,\ldots,d_\ell)$ if $D(\A,m) 
\simeq \oplus_{i=1}^\ell S[-d_i]$. 

(2)\,\,
We say that $(\A,m)$ is \textbf{tame}
if $\pd_S \Omega^p(\A,m) \le p$ for $0 \le p \le \ell$.
\label{multifree}
\end{define}

Let us gather several fundamental properties of these modules.

\begin{prop}[Saito's criterion, \cite{Z}]
Let $\theta_1,\ldots,\theta_\ell \in D(\A,m)$ be homogeneous derivations. Then they form a basis for $D(\A,m)$ if and only if 
for $M:=(\theta_i(x_j))_{1\le i,j \le \ell}$, it holds that 
$$
\det M=Q(\A,m)
$$
up to non-zero scalar. Equivalently, they form a basis if and only if $\theta_1,\ldots,\theta_\ell$ are independent over $S$ and $\sum_{i=1}^\ell \deg \theta_i=|m|$.
\label{saito}
\end{prop}

\begin{prop}[Lemmas 1.2 and 1.3, \cite{ATW}]
(1)\,\,If $(\A,m)$ is free, then $D^p(\A,m)=\wedge^p D(\A,m)$ for all $p$, and the basis for $D^p(\A,m)$ is obtained by taking the $p$-th wedges of the basis for $D(\A,m)$. 

(2)\,\, For all $\A$, it holds that 
$$
D^\ell(\A,m) \simeq S[-|m|].
$$
\label{wedge}
\end{prop}

\begin{prop}[Lemma 1.4, \cite{ATW}] 
Let $V=V_1 \times V_2$ and $\A_i$ be an arrangement in $V_i$. Let $m_i$ be a mulitplicity on $\A_i$, $m:=m_1\oplus m_2$ and let 
$\A=\A_1 \times \A_2$. Then 
$$
D^p(\A,m)=\oplus_{i+j=p} D^i(\A_1,m_1) \otimes_\K D^j(\A_2,m_2).
$$
\label{product}
\end{prop}

\begin{prop}[Proposition 2.8, \cite{A14}]
Let $(\A,m)$ be a multiarrangement. Then 
$$
D^p(\A,m) \simeq Q(\A,m) \Omega^{\ell-p}(\A,m)
$$
for $0 \le p \le \ell$. 
\label{identification}
\end{prop}

Now recall the Solomon-Terao bi-polynomial of a multiarrangement.

\begin{define}
For a multiarrangement $(\A,m)$, the \textbf{Solomon-Terao 
polynomial} $\Psi(\A,m;x,t)$ of $(\A,m)$ is defined by 
$$
\Psi(\A,m;x,t):=
\sum_{i=1}^\ell \mbox{Hilb}(D^p(\A,m);x))(t(x-1)-1)^p.
$$
Let $$
ST(\A,m;x):=\Psi(\A,m;x,-1)
$$
be the \textbf{Solomon-Terao polynomial} of $(\A,m)$.
\label{STmulti}
\end{define}

\begin{theorem}[Theorem 2.5, \cite{ATW}]
It holds that 
$$
\Psi(\A,m;x,t) \in \Q[x,t]
$$
and 
$$
ST(\A,m;x) \in \Q[x].
$$
Moreover, if $(\A,m)$ is free with $\exp(\A,m)=(d_1,\ldots,d_\ell)$, then 
$$
(-1)^\ell\Psi(\A,m;1,t)=\prod_{i=1}^\ell (t-d_i).
$$
\label{multiSTformula}
\end{theorem}

So we define 
$$
\chi(\A,m;t):=(-1)^\ell \Psi(\A,m;1,t)
$$
as a \textbf{characteristic polynomial} of $(\A,m)$.
To show Theorem \ref{multiSTformula}, the following was introduced.

\begin{define}[\cite{ATW}]
Let $\eta \in S_{d}$. Then define $\partial_p:D^p(\A,m) \rightarrow D^{p-1}(\A,m)$ by 
$$
\partial_p(\theta)(f_2,\ldots,f_p):=\theta(\eta,f_2,\ldots,f_p)
$$
for $f_2,\ldots,f_p \in S$ and $\theta \in D^p(\A,m)$. We call the complex 
$$
0 \rightarrow D^\ell(\A,m) \stackrel{\partial_\ell}{\rightarrow}
D^{\ell-1}(\A,m) \stackrel{\partial_{\ell-1}}{\rightarrow}
\cdots
\stackrel{\partial_2}{\rightarrow}
D(\A,m) \stackrel{\partial_1}{\rightarrow}S \rightarrow 0
$$
the \textbf{Solomon-Terao complex} of $(\A,m,\eta)$.
\label{STcomplex}
\end{define}

\begin{prop}[Lemma 2.4, \cite{ATW}]
Let $d>0$ be an integer. Then 
there is a non-empty Zariski open set $U_d =U_d(\A,m) \subset S_d$ such that all the 
homology group of the Solomon-Terao complex of $(\A,m,\eta)$ is of finite dimensional over $\K$ if $\eta \in U_d$.
\label{vanish}
\end{prop}

We say that $\A$ is \textbf{locally free along $H$} for $H \in \A$ if for any 
point $0\neq p \in H$, the localization
$$
\A_p:=\{L \in \A \mid p\in L\}
$$
is free. 
Let us introduce a result to check the tameness of deletions.

\begin{theorem}[Theorem 3.5, \cite{A14}]
Let $H \in \A$. 
If $\A$ is tame, locally free along $H$ and $\A^H$ is free, then $\A \setminus \{H\}$ is tame. In particular,
$\A \setminus \{H\}$ is tame if $\A$ and $\A^H$ are both free.
\label{tamedeletion}
\end{theorem}

For the example, let us introduce the following.

\begin{theorem}[Theorem 3.3, \cite{A5}]
Let $\A$ be free and $H \in \A$. Then there is a right exact sequence 
$$
0 \rightarrow \Omega^1(\A) 
\stackrel{\cdot \alpha_H}{\rightarrow}
\Omega^1(\A \setminus \{H\}) \rightarrow \Omega^1(\A^H) \rightarrow 0.
$$
\label{FST}
\end{theorem}

\section{Proofs}

Let us prove main results in this paper. First let us show Theorem \ref{regD}. For that we need the following, which is a multiarrangement version of what was introduced in \cite{AMMN}.

\begin{define}
Let $\eta \in U_d$ in the terminology of Proposition \ref{vanish}. Let 
$$
\mathfrak{a}(\A,m,\eta):=\{\theta(\eta) \in S \mid \theta \in D(\A,m)\}
$$
be the \textbf{Solomon-Terao ideal} of $(\A,m,\eta)$, and let 
$$
ST(\A,m,\eta):=S/\mathfrak{a}(\A,m,\eta)
$$
be the \textbf{Solomon-Terao algebra} of $(\A,m,\eta)$. 
\label{STalg}
\end{define}

As for the case of $m \equiv 1$, we can compute the Solomon-Terao polyomial by using the Solomon-Terao algebra when $\A$ is tame.

\begin{theorem}
Let $(\A,m)$ be tame and $\eta \in U_{d+1}=U_{d+1}(\A,m)$. Then 
$$
\Psi(\A,m;x,\frac{1-x^d}{x-1})=\mbox{Hilb}(\mbox{ST}(\A,m,\eta);x).
$$
In particular, when $d=1$, it holds that 
$$
ST(\A,m;x)=\mbox{Hilb}(ST(\A,m,\eta);x).
$$
\label{STequal}
\end{theorem}

\noindent
\textbf{Proof}.
This was proved in \cite{AMMN} when $m \equiv 1$, and the proof is essentially the same as it. We give a proof for the completeness. 
Namely, if $(\A,m)$ is tame, then the same argument as in \cite{AMMN} shows that 
$H_p(D^*(\A,m),\partial_*)=0$ unless $p=0$. So with Proposition \ref{vanish}, and by using the fact that $\partial_p$ 
is a map of degree plus $d$, we have 
\begin{eqnarray*}
\Psi(\A,m;x,\frac{1-x^d}{x-1})&=&\sum_{p=0}^\ell \mbox{Hilb}(D^p(\A,m);x)(-x^d)^p\\
&=&\sum_{p=0}^\ell \mbox{Hilb}(H_p(D^*(\A,m),\partial_*);x)(-x^d)^p\\
&=&\mbox{Hilb}(ST(\A,m,\eta);x)
\end{eqnarray*}
since $H_0(D^*(\A,m),\partial_*)=ST(\A,m,\eta)$. When $d=1$, the definition completes the proof of the second statement.
\owari
\medskip

Note that Theorem \ref{STequal} does not hold true in general if $(\A,m)$ is not tame, see Example \ref{ex2}. Since we will use the above new polynomial, let us give a name.

\begin{define}
Let $(\A,m)$ be tame and $\eta_{d+1} \in U_{d+1}(\A,m)$. Then denote 
$$
ST_{d+1}(\A,m;x):=\Psi(\A,m;x,\frac{1-x^d}{x-1})=\mbox{Hilb}(ST(\A,m,\eta_{d+1});x)
$$
which is called the \textbf{Solomon-Terao polynomial of order $d+1$}.
\label{STgeneral}
\end{define}

Note that $$
ST_2(\A,m;x)=ST(\A,m;x)
$$
and $ST_{d+1}(\A,m;x)$ is independent of the choice of a generic $\eta_{d+1} \in U_{d+1}$ by the above proof. When $d=1$ it was shown in \cite{AMMN} that 
$ST(\A;-1)=0$. Namely, 
$$
ST(\A;x)=(1+x)ST^+(\A;x).
$$
Let $ST^+(\A;x)$ be the \textbf{reduced Solomon-Terao polynomial} of $\A$.

Now recall the exact conjecture on the top degree monomials of the Solomon-Terao polynomials in \cite{AMMN}.

\begin{conj}[\cite{AMMN}, Conjecture 5.12]
Let $\eta \in U_{d+1}(\A)$. Then $ST_{d+1}(\A;x)$ is a monic polynomial of degree
$|\A|+\ell(d-1)$.
\label{AMMNconj}
\end{conj}

So our main result says that Conjecture \ref{AMMNconj} is true when $\A$ is tame. In fact, we prove the following more general result.

\begin{theorem}
Let $(\A,m)$ be tame and $\eta \in U_{d+1}(\A,m)$. Then $ST_{d+1}(\A,m;x)$ is a monic polynomial of degree
$|m|+\ell(d-1)$.
\label{AMMNmulti}
\end{theorem}

Let us start the proof of Theorem \ref{AMMNmulti}. For that, first let us complete the proof of the regularity.
\medskip

\noindent
\textbf{Proof of Theorem \ref{regD}}.
We prove by applying the argument as in \cite{MS} with a slight modification for multiarrangements. We prove by induction on $\ell \ge 1$. When $\ell=1$ there is nothing to show. 
Let $\ell=2$. Since $D(\A,m)$ is reflexive, $D(\A,m)$ is free. So for $D(\A,m)$, it suffices to show that 
the larger degree of the basis for $D(\A,m)$ is at most $n-1$, where $n:=|m|$. If not, then by Proposition \ref{saito}, 
the smaller degree is $0$, which is impossibile since $\A$ is essential. Since $D^2(\A,m) \simeq S[-n]$, the case when 
$\ell=2$ is completed.

Assume that the statement is true up to $\ell-1 \ge 2$ and consider a multiarrangement $(\A,m)$ in $\K^\ell$. 
We prove by induction on $|m|$. If $|m|\le \ell$, then $m \equiv 1$ since $\A$ is essential, and $\A$ is free with $\exp(\A)=(1,\ldots,1)$. So 
Theorem \ref{wedge} completes the proof. Assume that the statement is true up to $|m|<n $ and show the case when $|m|=n$. 
Let $H_1,\ldots,H_\ell \in \A$ be independent hyperplanes and let $\delta_i$ be a multiplicity on $\A$ such that $\delta_i(H)=1$ only 
when $H=H_i$, and zero otherwise. Let $m_i:=m-\delta_i$ and let $m':=m-\sum_{i=1}^\ell \delta_i$. If $m(H_i)>1$ for all $i$, then 
the deleted $(\A,m_i)$ are still all essential, and 
$$
\alpha_{H_i} D^p(\A,m_i) \subset D^p(\A,m) \subset D^p(\A,m_i)
$$
for $i=1,\ldots,\ell$. Since the induction hypothesis shows that 
$\reg D^p(\A,m_i) \le |m|-1-\ell+p$, Theorem \ref{DS} completes the proof. Assume that there is $i$, say $i=1$ such that $m(H_1)=1$. Then 
for $\A_1:=\A \setminus \{H_1\}$, it may hold that $\A_1$ is not essential. Let $H_1:x_1=0$. Then this means that $Q(\A)=x_1Q(\A_1)$ with 
$Q(\A_1) \in \K[x_2,\ldots,x_\ell]$. So Proposition \ref{product} shows that 
$$
D^p(\A,m)= S D^p(\A_1^e,m_1) \oplus S  D^{p-1}(\A_1^e,m_1) \wedge x_1 \partial_{x_1}.
$$
Here $\A_1^e$ is an arrangement in $x_1=0$ obtained by dividing $\A_1$ by $x_2=\cdots=x_\ell=0$.
Since $\reg D^p(\A_1^e,m_1) \le n-1-(\ell-1)+p=n-\ell+p$ and 
$\reg D^{p-1}(\A_1^e,m_1) \wedge x_1 \partial_{x_1} \le (n-1-(\ell-1)+p-1)+1
=n-\ell+p$ by induction, it holds that $\reg D^p(\A,m) \le n-\ell+p$.

Let us show the statement for the logarithmic differential $p$-forms. 
Note that $Q(\A,m)\Omega^p(\A,m) \simeq D^{\ell-p}(\A,m)$ by Proposition \ref{identification}. Thus 
$$
\reg \Omega^p(\A,m) =\reg D^{\ell-p}(\A)[-n] \le n-\ell+\ell-p-n=-p,
$$
which completes the proof. 

Note that we can show $\reg \Omega^p(\A,m) \le -p$ by the same way as \cite{MS} by using Theorem \ref{DS}, which gives another proof of the result in \cite{B}, as follows. Let $M:=Q(\A,m) \Omega^p(\A,m)$ and let $M_i:=Q(\A,m-\delta_{H_i})\Omega(\A,m-\delta_{H_i})$. 
Since $\alpha_{H_i} \Omega^p(\A,m) \subset \Omega^p(\A,m-\delta_{H_i})$, we know that 
$$
\alpha_{H_i} M_i \subset M \subset M_i
$$
for $i=1,\ldots,\ell$. Thus we can apply the same argument as above, and combining $\reg Q(\A,m) \Omega^p(\A,m)=|m|+\reg \Omega^p(\A,m)$, we complete the proof.
\owari
\medskip

By using this regularity bound, we investigate the coefficients of the first and second smallest/largest powers of $x$ in $ST_{d+1}(\A,m;x)$. 

\begin{prop}
Assume that $(\A,m)$ is tame and irreducible. Then 
$$
ST_{d+1}(\A,m;x)=
\sum_{i=0}^d \displaystyle \frac{(\ell+i-1)!}{(\ell-1)!} x^i
+
(\displaystyle \frac{(\ell+d)!}{(\ell-1)!}-d_2)x^{d+1}
+\sum_{i=d+2}^\infty c_p x^p.
$$
Here $d_2=1$ if $m \not \equiv 1$ and $d_2=0$ otherwise.
\label{const}.
\end{prop}

\noindent
\textbf{Proof}.
We use Theorem \ref{STalg}, i.e., we compute by using the Solomon-Terao algebra. 
Note that $\partial_p$ is a map of degree plus $d$. So $\dim_\K
ST(\A,m,\eta_{d+1})_p=\dim_\K S_p$ for $0 \le p \le d$ since $D(\A,m)_{\le 0}=(0)$ by the irreducibility. So the coefficients of $ST_{d+1}(\A,m;x)$ of $x^i, 0 \le i \le d$, are clear.

Next consider the coefficient of $x^{d+1}$. By the same reason, $\dim_\K S_{d+1}=\displaystyle \frac{(\ell+d)!}{(\ell-1)!}$. 
So it suffices to compute $\dim_\K \mathfrak{a}(\A,m,\eta_{d+1})_{d+1}$ which comes from $D(\A,m)_1$.
Since $\A$ is essential, $\dim_\K D(\A)_1=1$ since there is the Euler derivation at degree one. If there is $H \in \A$ such that $m(H) \ge 2$, then the Euler derivation cannot be tangent to any hyperplane with multiplicity 
more than or equal to two. So $D(\A,m)_1=(0)$, which completes the proof.\owari
\medskip

Now let us show more general theorem than Theorem \ref{main1}.

\begin{theorem}
Let $\A$ be essential and $(\A,m)$ be tame. Then 
$
ST_{d+1}(\A,m;x)
$ is a monic polynomial of degree $r:=|m|+\ell(d-1)$.
\label{mainmulti}
\end{theorem}

\noindent
\textbf{Proof}.
Let $n:=|m|$. 
Recall that 
$$
\Psi(\A,m;x,t)=
\sum_{p=0}^\ell 
\mbox{Hilb}(D^p(\A,m);x)(t(x-1)-1)^p.
$$
Let $\mbox{Hilb}(D^p(\A,m);x)=\displaystyle \frac{f_p(x)}{(1-x)^\ell}$.
Since $D^\ell(\A,m)\simeq S[-n]$, it holds that $f_\ell(x)=x^n$.
By the tameness of $(\A,m)$, it holds that $\pd_S D^p(\A,m) \le \ell-2$. So by 
Theorem \ref{regD} and Proposition \ref{identification}, we know that $\deg f_p(x) \le \ell-p+n-\ell+p=n$ for all $p$. Hence 
$$
f_p(x)=a_px^n+(\mbox{lower terms}),
$$
thus
\begin{eqnarray}
ST_{d+1}(\A,m;x)&=&\Psi(\A,m;x,\frac{1-x^d}{x-1})=
\sum_{p=0}^\ell \displaystyle \frac{f_p(x)}{(1-x)^\ell}(-x^d)^p\nonumber \\
&=&\displaystyle \frac{x^{n+d\ell}-a_{\ell-1}x^{n+d(\ell-1)}+
(\mbox{lower degree terms})
}{(x-1)^\ell}.
\label{eq1}
\end{eqnarray}
Since this is a polynomial, we know that the top degree term is $x^{n+\ell(d-1)}=x^r$, which completes the proof.\owari
\medskip

\noindent
\textbf{Proof of Theorem \ref{main1}}. 
Apply Theorem \ref{mainmulti} to the case when $m \equiv 1$
and $d=1$. \owari
\medskip

\noindent
\textbf{Proof of Corollary \ref{3topdegree}}.
Apply Theorem \ref{main1} to three dimensional arrangements 
when and $d=1$. \owari
\medskip

Now let us consider the second highest coefficient, whose behavior differs depending on $d=1$ or not. First consider when $d=1$, so the second highest degree is $|m|-1$. Let us show the following.

\begin{theorem}
Let $(\A,m)$ be irreducible and tame.
Then for $n:=|m|$, the coefficient of $x^{n-1}$ in $ST(\A,m;x)$ is 
$\ell+a$ with $a \ge 0$. Here $a$ is the number of relations of degree 
$n$ among a minimal set of generators for $D^{\ell-1}(\A,m)$, which is by Proposition \ref{identification}, equivalent to the relations of degree $0$ for $\Omega^1(\A,m)$.
\label{secondhighest}
\end{theorem}

\noindent
\textbf{Proof}.
Let us compute $\Psi(\A,m;x,t)$ directly. By the proof of Theorem \ref{mainmulti}, we can express  $f_p(x)=a_px^n+b_px^{n-1}+\cdots$ in that terminology. Then 
we get 
\begin{eqnarray*}
\Psi(\A,m;x,t)&=&
(-1)^\ell t^\ell x^n+(-1)^{\ell-1}t^{\ell-1}
\displaystyle\frac{\ell-a_{\ell-1}}{x-1}x^n\\
&+&
(-1)^{\ell-1}t^{\ell-1}
\displaystyle\frac{-b_{\ell-1}}{x-1}x^{n-1}
+(-1)^{\ell-2}t^{\ell-2}
\displaystyle\frac{\frac{\ell(\ell-1)}{2}-(\ell-1) a_{\ell-1}+a_{\ell-2}}{(x-1)^2}x^n\\
&=&
(-1)^\ell t^\ell x^n+(-1)^{\ell-1}t^{\ell-1}
\displaystyle\frac{(\ell-a_{\ell-1})x^n-b_{\ell-1}x^{n-1}}{x-1}\\
&+&
(-1)^{\ell-2}t^{\ell-2}
\displaystyle\frac{\frac{\ell(\ell-1)}{2}-(\ell-1) a_{\ell-1}+a_{\ell-2}}{(x-1)^2}x^n
\end{eqnarray*}
by ignoring the terms $t^i$ for $i \le \ell-3$. So 
\begin{eqnarray*}
ST(\A,m;x)&=&
x^n+(\ell-a_{\ell-1})x^{n-1}+
(\frac{\ell(\ell+1)}{2}-\ell a_{\ell-1}- b_{\ell-1}+a_{\ell-2})x^{n-2}\\
&+&
\mbox{(lower degree terms)}.
\end{eqnarray*}
So it suffices to show that $a:=a_{\ell-1} \le 0$. Note that $a$ comes from degree 
$n$-part of $\mbox{Hilb}(D^{\ell-1}(\A,m);x)$, $\pd_S D^{\ell-1}(\A,m) \le 1$ and 
$\reg D^{\ell-1}(\A,m) \le n-1$ by the tameness and Theorem \ref{regD}. So 
for $\mbox{Hilb}(D^{\ell-1}(\A,m);x)=a_{\ell-1}x^n+\cdots$, $a_{\ell-1}$ is nothing but the absolute value of the number of 
degree-$n$ relations among a minimal set of generators for $D^{\ell-1}(\A,m)$. Since a minimal 
free resolution of $D^{\ell-1}(\A,m)$ is of the form
$$
0 \rightarrow 
\oplus_{k \ge -n+1} S[-k]^{s_k} \oplus S[-n]^{-a_{\ell-1}}  \rightarrow F_0 \rightarrow D^{\ell-1}(\A,m) \rightarrow 0,
$$
we complete the proof by combining Proposition \ref{identification}.\owari
\medskip

\noindent
\textbf{Proof of Theorem \ref{main2}}.
Apply Theorem \ref{secondhighest} to the case $m \equiv 1$. \owari
\medskip

The proof of Theorem \ref{secondhighest} implies the following.

\begin{cor}
\begin{eqnarray*}
\Psi(\A,m;x,t)&=&(-1)^\ell t^\ell x^n+
(-1)^{\ell-1}t^{\ell-1}((\ell-a_{\ell-1})x^{n-1}-b_{\ell-1}x^{n-2})
\end{eqnarray*}
by ignoring the terms $t^i$ for $i \le \ell-3$.
\label{tl}
\end{cor}

\begin{rem}
Note that if $(\A,m)$ is free, then $a_p=0$ since all $D^p(\A,m)$ are free by Proposition \ref{wedge}, thus no relations. 
However, if $m \equiv 1$ and for a basis $\theta_1=\theta_E,\theta_2,\ldots,\theta_\ell$ for $D(\A)$, there is a 
derivation 
$\theta_2 \wedge \cdots \wedge \theta_\ell$
in the basis for $D^{\ell-1}(\A)=\wedge^{\ell-1}D(\A)$ by Proposition \ref{wedge} which is of degree $|\A|-1$. So in this case 
$b_{\ell-1}=1$ and $ST(\A;x)$ is palindromic for the first/last three terms.
\end{rem}

\begin{rem}
In the statement of Theorem \ref{secondhighest}, the tameness for $\A$ is necessary. In fact, as in Example \ref{ex2}, there is a non-tame arrangement 
such that the coefficient of $x^{|\A|-1}$ is strictly smaller than $\ell$.
\end{rem}

Next let us consider when $d>1$. In this case we have the following.

\begin{prop}
Let $\eta \in U_{d+1}(\A,m)$ with $d>1$ and $(\A,m)$ be tame. Let $n:=|m|$ and 
$r:=n+\ell(d-1)$. Then 
\begin{eqnarray*}
ST_{d+1}(\A,m;x)&=&\sum_{i=0}^d \displaystyle \frac{(\ell+i-1)!}{(\ell-1)!} x^i
+
(\displaystyle \frac{(\ell+d)!}{(\ell-1)!}-d_2)x^{d+1}+
\cdots\\
&+&
(\displaystyle \frac{(\ell+d-1)!}{(\ell-1)!}-a_{\ell-1})x^{r-d}+
\sum_{i=0}^{d-1}\displaystyle \frac{(\ell+i-1)!}{(\ell-1)!}x^{r-i}.
\end{eqnarray*}
Here $a_{\ell-1}$ and so on are the same as in the proof of Theorem \ref{mainmulti} and equation \ref{eq1}.
\end{prop}

\noindent
\textbf{Proof}. Clear by the equation (\ref{eq1}). \owari
\medskip

Since all three dimensional arrangements are tame, we can show the following.

\begin{cor}
Let $\ell=3$. Then
\begin{eqnarray*}
ST_{d+1}(\A,m;x)&=&\sum_{i=0}^d \displaystyle \frac{(\ell+i-1)!}{(\ell-1)!} x^i
+
(\displaystyle \frac{(\ell+d)!}{(\ell-1)!}-d_2)x^{d+1}+
\cdots\\
&+&
(\displaystyle \frac{(\ell+d-1)!}{(\ell-1)!}-a_{\ell-1})x^{|m|+d-2}+
\sum_{i=0}^{d-1}\displaystyle \frac{(\ell+i-1)!}{(\ell-1)!}x^{|m|+2d-2-i}.
\end{eqnarray*}
\end{cor}

\section{Examples}

In this section let us see/use what we proved to Solomon-Terao polynomials 
of several arrangements.

\begin{example}
Let $\A$ be an irreducible arrangement in $\R^3$ defined by 
$$
xyz(x+y+z)=0.
$$
Then 
$$
\Psi(\A;x,t)=-x^4t^3+4x^3t^2-t(5x^2+x)+2x+1,
$$
see \cite{A12}, Example 6.5 for example. 
Hence $$
ST(\A;x)=x^4+4x^3+5x^2+3x+1.
$$
Thus Theorem \ref{secondhighest} implies that there is one relation for the minimal set of generators for 
$D^2(\A) \simeq Q(\A) \Omega^1(\A)$ at degree $4=|\A|$.
In fact, $\Omega^1(\A)$ has a minimal free resolution 
$$
0 \rightarrow 
S[0] \rightarrow  S[1]^{4} \rightarrow \Omega^1(\A) \rightarrow 0,
$$
see \cite{RT} or \cite{A5}. Thus Proposition \ref{identification} shows that 
$$
0 \rightarrow S[-4] \rightarrow  S[-3]^{4} \rightarrow D^2(\A) \rightarrow 0,
$$
as predicted in Theorem \ref{secondhighest}.
\label{ex1}
\end{example}

\begin{example}
Let $\A$ be an irreducible generic arrangement in $\K^\ell$ with $|\A|=n$, here $\K$ is a field of characteristic zero. Then it was proved in Theorem 3.2.1, \cite{RT} that $\Omega^1(\A)$ has a minimal free resolution 
$$
0 \rightarrow S[0]^{n-\ell} \rightarrow S[1]^n \rightarrow \Omega^1(\A) \rightarrow 0.
$$
Since $Q(\A)\Omega^1(\A) \simeq D^{\ell-1}(\A)$ as in Proposition \ref{identification}, we have a minimal free resolution
$$
0 \rightarrow S[-n]^{n-\ell} \rightarrow S[-n+1]^n \rightarrow \Omega^1(\A) \rightarrow 0.
$$
So Theorem \ref{secondhighest} shows that, since $\A$ is tame, 
$$
ST(\A;x)=1+\ell x+\cdots+n x^{n-1}+x^n.$$
\end{example}

\begin{example}
Let $\A$ be a free arrangement in $\R^4$ defined by 
$$
Q(\A)=x_1x_2x_3x_4(x_1+x_2)(x_1+x_3)(x_1+x_4)(x_1+x_2+x_3)(x_1+x_2+x_4)(x_1+x_3+x_4)=0.
$$
It is known that $\exp(\A)=(1,3,3,3)$, see \cite{A5}, Example 6.1 for example.

First let $H:x_2=0$ and $\A':=\A \setminus \{H\}$.
Then Theorem \ref{tamedeletion} shows that $\A'$ is tame since $\A$ are both $\A^H$ free. So we can apply 
Theorem \ref{secondhighest} to $\A'$. Then we can compute that 
$$
ST(\A';x)=1+4x+9x^2+16x^3+21x^4+21x^5+17x^6+10x^7+4x^8+x^9.
$$
Hence Theorem \ref{secondhighest} tells us that there are no degree-$9$ relations among minimal set of generators for 
$D^{\ell-1}(\A') \simeq Q(\A') \Omega^1(\A')$.

Contrary to it, let $L:x_1=0$ and consider $\B:=\A \setminus \{L\}$.
Then we can compute that 
$$
ST(\B;x)=1+4x+9x^2+16x^3+21x^4+21x^5+17x^6+9x^7+3x^8+x^9.
$$
So the coefficient of $x$ is strictly larger than that of $x^8$, i.e., the statement in Theorem \ref{secondhighest} does not hold true. 
This means that $\B$ is not tame, i.e., $\pd_S \Omega^1(\B)=\pd_S D^{\ell-1}(\B)=2>1$.
We can check that $\B$ is not tame by a different way. By Theorem \ref{FST}, it holds that 
$$
0 \rightarrow 
\Omega^1(\A) \stackrel{\cdot \alpha_L}{\rightarrow }
\Omega^1(\B) \rightarrow \Omega^1(\A^L) \rightarrow 0.
$$
As in Example 6.2 in \cite{A5}, we know that $\pd_S \Omega^1(\A^L)=2$. Since $\pd_S \Omega^1(\A)=0$, the Ext-long exact 
sequence shows that $\pd_S \Omega^1(\B)=2$. 
As a consequence, $\B$ is not tame.
\label{ex2}
\end{example}

\end{document}